# Even and Odd Pairs of Lattice Paths With Multiple Intersections

Ira M. Gessel[*]
Department of Mathematics
Brandeis University
Waltham, MA 02254-9110

and

Walter Shur
11 Middle Road
Port Washington, NY 10010

**Abstract.** Let $M_{r,s}^{n,k}$ be the number of ordered pairs of paths in the plane, with unit steps E or N, that intersect k times in which the first path ends at the point (r,n-r) and the second path ends at the point (s,n-s). Let

$$N_E(n,k,p) = \sum_{r+s=2p} M_{r,s}^{n,k}$$

and

$$N_O(n,k,p) = \sum_{r+s=2p+1} M_{r,s}^{n,k}.$$

We study the numbers $M_{r,s}^{n,k}$, $N_k^{n,r} = M_{r,r}^{n,k}$, $N_E(n,k,p)$, and $N_O(n,k,p)$, prove several simple relations among them, and derive a simpler formula for $M_{r,s}^{n,k}$ than appears in [1].

Introduction

A recent paper [1] considered paths that begin at the SW corner of a lattice rectangle and proceed with unit steps in either of the directions E or N. Let $M_{r,s}^{n,k}$ be the number of ordered pairs of such paths that intersect k times in which the first path ends at the point (r,n-r) and the second path ends at the point (s,n-s), and let $N_k^{n,r} = M_{r,r}^{n,k}$. Intersections at the origin and endpoints of the paths are not counted. It is clear that $M_{r,s}^{n,k} = M_{s,r}^{n,k}$ and that $N_{n-1}^{n,r} = \binom{n}{r}$ for $0 \le r \le n$.

---

[*] partially supported by NSF grant DMS-9306297

In [1] the following formulas for $M_{r,s}^{n,k}$ and $N_k^{n,r}$ are given:

$$M_{r,s}^{n,k} = 2\sum_t \sum_j (-1)^j \frac{(s-j-r+1+2t)}{n-1-j-2t} \binom{k}{2t+1}\binom{k-1-2t}{j}\binom{n-1-j-2t}{s-j}\binom{n-1-j-2t}{r-1-2t}$$

$$+ \frac{s-r}{n-k} \sum_j \binom{k}{j}\binom{n-k}{r-j}\binom{n-k}{s-j}, \quad (0 \le k \le n-1, \; r < s)$$

$$N_k^{n,r} = \frac{2(k+1)}{n-k-1} \sum_i \binom{k}{i}\binom{n-k+i-1}{r}\binom{n-i-1}{n-r}, \quad (0 \le k \le n-2)$$

$$N_k^{n,r} = \frac{2(k+1)}{r} \sum_i (-1)^i \frac{\binom{k}{i}\binom{k-i}{i}\binom{n-i-2}{r-1}\binom{n-i-1}{r-i-1}}{\binom{n-i-2}{i}}, \quad (0 \le k \le n-2)$$

Our main object of study in this paper is the sum of the numbers $M_{r,s}^{n,k}$ over r and s where r+s is fixed. We are particularly interested in considering even and odd values of r+s separately. Thus we define

$$N_E(n,k,p) = \sum_{r+s=2p} M_{r,s}^{n,k}$$

and

$$N_O(n,k,p) = \sum_{r+s=2p+1} M_{r,s}^{n,k}.$$

We will prove the following three formulas for these numbers:

**Theorem 1**

$$N_E(n,k,p) = \frac{n}{k+1} N_k^{n,p}$$

**Theorem 2**

$$N_O(n,k,p) = 2\left[\binom{n-1}{p}^2 - \sum_{i=0}^{k-2} N_i^{n-1,p}\right], \quad n > 0$$

**Theorem 3**

$$\sum_p N_O(n,k,p) = \sum_p N_E(n,k,p) = 2^{k+1}\binom{2n-k-2}{n-1}$$

We will also prove a recurrence for $N_k^{n,r}$ and a simpler formula for $M_{r,s}^{n,k}$:

**Theorem 4**

$$\frac{n-k-1}{k+1} N_k^{n,r} = \frac{n-2k}{k} N_{k-1}^{n,r} + N_{k-2}^{n-1,r} + N_{k-2}^{n-1,r-1}, \quad k \geq 1$$

**Theorem 5**

$$M_{r,s}^{n,k} = \sum_{i=0}^{k} \binom{k}{i}\left[\binom{n-k+i-1}{s-1}\binom{n-i-1}{n-r-1} - \binom{n-k+i-1}{s}\binom{n-i-1}{n-r}\right], \quad r<s$$

**The Generating Function for $M_{r,s}^{n,k}$**

To find a generating function for the numbers $M_{r,s}^{n,k}$, we use a decomposition of pairs of paths counted by $M_{r,s}^{n,k}$ for r<s into s-r+1 parts. To describe this decomposition, let $d_m$, for m=0,1,...,n, be the (nonnegative) difference in x-coordinates of the mth points on the two paths. For each m>0, $d_m - d_{m-1} = 0$ or 1. Since $d_n = s-r$ and $s-r \leq n$, for each j=0,1,...,s-r there is at least one value of m such that $d_m = j$. Let $m_j$ be the largest such value of m. Then, if we let $m_{-1} = 0$, the jth part of our decomposition, for j=0,1,...,s-r, consists of the steps of the two paths from the $m_{j-1}$th points to the $m_j$th points (see Figure 1).

Note that the zeroth part of the decomposition will be empty if $m_0 = 0$. If it is not empty, it is a pair of paths counted by $N_{k-1}^{n',r'}$ for some $n'$ and $r'$. Each other part of the decomposition can be transformed, by translating its paths so that they start at the origin, to a pair of paths counted by $M_{r',r'+1}^{n',0}$ for some $n'$ and $r'$.

**Figure 1**

Moreover, the number of intersections in the pair of paths is zero if the zeroth part is empty, and is one more than the number of intersections in the zeroth part if the zeroth part is nonempty.

The decomposition gives a factorization for the generating function for $M_{r,r+j}^{n,k}$. We assign the weight $x_1^r x_2^{r+j} y_1^{n-r} y_2^{n-r-j}$ to an ordered pair of paths, each path starting at the origin, the first ending at the point $(r, n-r)$ and the second ending at the point $(r+j, n-r-j)$.

Then for fixed $j>0$, we have

$$\sum_{\substack{n \geq j \\ 0 \leq r \leq n-j}} M_{r,r+j}^{n,0} x_1^r x_2^{r+j} y_1^{n-r} y_2^{n-r-j} = \left[ \sum_{\substack{n \geq 1 \\ 0 \leq r < n}} M_{r,r+1}^{n,0} x_1^r x_2^{r+1} y_1^{n-r} y_2^{n-r-1} \right]^j$$

$$\sum_{\substack{n \geq j \\ 0 \leq r \leq n-j}} M_{r,r+j}^{n,k} x_1^r x_2^{r+j} y_1^{n-r} y_2^{n-r-j} = \left[ \sum_{\substack{n \geq 1 \\ 0 \leq r \leq n}} N_{k-1}^{n,r} x_1^r x_2^r y_1^{n-r} y_2^{n-r} \right] \left[ \sum_{\substack{n \geq 1 \\ 0 \leq r < n}} M_{r,r+1}^{n,0} x_1^r x_2^{r+1} y_1^{n-r} y_2^{n-r-1} \right]^j , \quad k>0.$$

It follows from [1] that for fixed $k \geq 1$ we have

$$\sum_{\substack{n \geq 1 \\ 0 \leq r \leq n}} N_{k-1}^{n,r} x_1^r x_2^r y_1^{n-r} y_2^{n-r} = [x_1 x_2 + y_1 y_2 + 2f(x_1 x_2, y_1 y_2)]^k \tag{1}$$

where $f(x,y)$ satisfies $f=(x+f)(y+f)$ and $f(0,0)=0$. It is clear that $N_0^{n,r} = 2 M_{r-1,r}^{n-1,0}$, and thus from the case $k=1$ of (1) we get

$$\sum_{\substack{n \geq 1 \\ 0 \leq r < n}} M_{r,r+1}^{n,0} x_1^r x_2^{r+1} y_1^{n-r} y_2^{n-r-1} = \frac{f(x_1 x_2, y_1 y_2)}{x_1 y_2} .$$

Therefore, for $k \geq 0$,

$$\sum_{\substack{n \geq j \\ 0 \leq r \leq n-j}} M_{r,r+j}^{n,k} x_1^r x_2^{r+j} y_1^{n-r} y_2^{n-r-j} = [x_1 x_2 + y_1 y_2 + 2f(x_1 x_2, y_1 y_2)]^k \left[ \frac{f(x_1 x_2, y_1 y_2)}{x_1 y_2} \right]^j . \tag{2}$$

Summing over all $j \geq 1$, we obtain

$$\sum_{\substack{n \geq 1 \\ r < s}} M_{r,s}^{n,k} x_1^r x_2^s y_1^{n-r} y_2^{n-s} = [x_1 x_2 + y_1 y_2 + 2f(x_1 x_2, y_1 y_2)]^k \left[ \frac{\frac{f(x_1 x_2, y_1 y_2)}{x_1 y_2}}{1 - \frac{f(x_1 x_2, y_1 y_2)}{x_1 y_2}} \right] . \tag{3}$$

## The Generating Functions for $N_E(n,k,p)$ and $N_O(n,k,p)$.

Setting $x_1 = x_2 = x$ and $y_1 = y_2 = y$ in (3) yields

$$\sum_{\substack{n \geq 1 \\ t \geq 1}} \left( \sum_{\substack{r+s=t \\ r<s}} M_{r,s}^{n,k} \right) x^t y^{2n-t} = (x^2 + y^2 + 2f(x^2, y^2))^k \frac{\frac{f(x^2, y^2)}{xy}}{1 - \frac{f(x^2, y^2)}{xy}} .$$

Noting that

$$\frac{\frac{f(x^2, y^2)}{xy}}{1 - \frac{f(x^2, y^2)}{xy}} = \frac{\left[\frac{f(x^2, y^2)}{xy}\right]^2}{1 - (\frac{f(x^2, y^2)}{xy})^2} + \frac{\frac{f(x^2, y^2)}{xy}}{1 - (\frac{f(x^2, y^2)}{xy})^2} ,$$

we readily obtain:

$$\sum_{\substack{n \geq 1 \\ p \geq 1}} \left( \sum_{\substack{r+s=2p \\ r<s}} M_{r,s}^{n,k} \right) x^p y^{n-p} = (x+y+2f(x,y))^k \frac{\frac{f(x,y)^2}{xy}}{1 - \frac{f(x,y)^2}{xy}}$$

and

$$\sum_{\substack{n \geq 1 \\ p \geq 1}} \left( \sum_{\substack{r+s=2p+1 \\ r<s}} M_{r,s}^{n,k} \right) x^p y^{n-p} = (x+y+2f(x,y))^k \frac{\frac{f(x,y)}{x}}{1 - \frac{f(x,y)^2}{xy}} .$$

By definition,

$$N_E(n,k,p) = \sum_{r+s=2p} M_{r,s}^{n,k} = 2 \sum_{\substack{r+s=2p \\ r<s}} M_{r,s}^{n,k} + N_k^{n,p}$$

and

$$N_O(n,k,p) = \sum_{r+s=2p+1} M_{r,s}^{n,k} = 2 \sum_{\substack{r+s=2p+1 \\ r<s}} M_{r,s}^{n,k} .$$

Hence, writing f for f(x,y), and using (1) to obtain the generating function for $N_k^{n,p}$, we find the generating functions

$$\sum_{n,p} N_E(n,k,p) x^p y^{n-p} = 2(x+y+2f)^k \frac{\frac{f^2}{xy}}{1-\frac{f^2}{xy}} + (x+y+2f)^{k+1} \tag{4}$$

and

$$\sum_{n,p} N_O(n,k,p) x^p y^{n-p} = 2(x+y+2f)^k \frac{\frac{f}{x}}{1-\frac{f^2}{xy}} . \tag{5}$$

**Proof of Theorem 1**

We need to show that the generating function for the numbers $\frac{n}{k+1} N_k^{n,p}$ is equal to (4). From (1), we know that $\sum_{n,p} N_k^{n,p} x^p y^{n-p} = [x+y+2f]^{k+1}$ and that $f=(x+f)(y+f)$. Noting that $(x\frac{\partial}{\partial x}+y\frac{\partial}{\partial y})x^p y^{n-p}=nx^p y^{n-p}$, and deriving the facts that

$\frac{\partial f}{\partial x}=\frac{y+f}{1-x-y-2f}$ and $\frac{\partial f}{\partial y}=\frac{x+f}{1-x-y-2f}$, we have

$$\sum_{n,p} \frac{n}{k+1} N_k^{n,p} x^p y^{n-p} = \frac{1}{k+1} (x\frac{\partial}{\partial x}+y\frac{\partial}{\partial y})(x+y+2f)^{k+1}$$

$$= (x+y+2f)^k \frac{2xy+x+y-x^2-y^2}{1-x-y-2f} .$$

Using the relationships $f^2=f-xf-yf-xy$ and $\frac{f}{xy-f^2}=\frac{1}{1-x-y-2f}$, we see that this is equal to (4).

**Proof of Theorem 2**

Theorem 2 rearranged states that

$$\frac{1}{2} N_O(n,k,p) + \sum_{i=0}^{k-2} N_i^{n-1,p} = \binom{n-1}{p}^2 ,$$ where $n \geq 1$ and the sum is taken as 0 if k<2. We need to show that the generating function for the numbers on the left side of this equation equals $\sum_{\substack{n,p \\ n \geq 1}} \binom{n-1}{p}^2 x^p y^{n-p}$.

Adding the appropriate generating functions, we have,

$$G_k(x,y) = \sum_{n,p} [\frac{1}{2}N_o(n,k,p) + \sum_{i=0}^{k-2} N_i^{n-1,p}] x^p y^{n-p} = (x+y+2f)^k \frac{yf}{xy-f^2} + \sum_{i=0}^{k-2} y(x+y+2f)^{i+1}.$$

Define $\Delta_k$ by $\Delta_k h(k) = h(k+1) - h(k)$. We have

$$\Delta_k G_k(x,y) = -y, \quad k=0,$$

$$\Delta_k G_k(x,y) = \frac{yf}{xy-f^2}(x+y+2f)^k(x+y+2f-1) + y(x+y+2f)^k, \quad k \geq 1.$$

Since $\frac{f}{xy-f^2} = \frac{1}{1-x-y-2f}$, we have for $k \geq 1$, $\Delta_k G_k(x,y) = 0$ and

$$G_k(x,y) = G_1(x,y).$$

Since $G_0(x,y) = G_1(x,y) + y$, all we need to do is deal with $G_0(x,y)$ and show that

$$\frac{yf}{xy-f^2} = \sum_{n,p} \binom{n-1}{p}^2 x^p y^{n-p}.$$

(Note that the term $y$ in $G_0(x,y) = G_1(x,y) + y$ counts the number $\frac{1}{2}N_O(1,0,0) = \binom{0}{0}^2$.)

From [1], we have

$$u_0(y, \frac{x}{y}) = \sum_{n,r} N_k^{n,r} x^r y^{n-r} = x+y+2f(x,y),$$

$$\frac{1}{1-u_0(y,x/y)} = \sum_{n,p} \binom{n}{p}^2 x^p y^{n-p}.$$

Hence, since $\frac{yf}{xy-f^2} = \frac{y}{1-x-y-2f}$, we have

$$\frac{yf}{xy-f^2} = \sum_{n,p} \binom{n-1}{p}^2 x^p y^{n-p}.$$

**Proof of Theorem 3**

We first show that the generating functions for $\sum_p N_E(n,k,p)$ and $\sum_p N_O(n,k,p)$ are equal. We obtain them by letting $y=x$ in (4) and (5).

Thus, we have

$$\sum_n [\sum_p N_E(n,k,p)] x^n = (2x+2f)^k \frac{2x}{1-2x-2f}$$

$$\sum_n [\sum_p N_O(n,k,p)] x^n = 2(2x+2f)^k \frac{xf}{x^2-f^2} = (2x+2f)^k \frac{2x}{1-2x-2f}.$$

From Theorem 1, we have $\sum_p N_E(n,k,p) = \frac{n}{k+1} \sum_p N_k^{n,p}$. From [1] we know that $\sum_p N_k^{n,p} = 2^{k+1}(k+1) \frac{(2n-k-2)!}{n!(n-k-1)!}$. Hence,

$$\sum_p N_E(n,k,p) = \sum_p N_O(n,k,p) = 2^{k+1} \binom{2n-k-2}{n-1}.$$

**Proof of Theorem 4**

We begin with the following recurrences, where $k<n-1$,

$$M_{r,s}^{n,k} = M_{r,s}^{n-1,k} + M_{r-1,s}^{n-1,k} + M_{r,s-1}^{n-1,k} + M_{r-1,s-1}^{n-1,k}, \quad r<s-1$$

$$M_{r,r+1}^{n,k} = M_{r,r+1}^{n-1,k} + M_{r-1,r+1}^{n-1,k} + N_{k-1}^{n-1,r} + M_{r-1,r}^{n-1,k}, \quad r=s-1.$$

Apply $2\sum_{\substack{r+s=2p+1 \\ r<s}}$ to both sides of the first recurrence, multiply both sides of the second recurrence by 2 and set $r=p$, and add the resulting recurrences together to obtain

$$N_O(n,k,p) = N_O(n-1,k,p) + [N_E(n-1,k,p) - N_k^{n-1,p}]$$
$$+ [N_E(n-1,k,p) + 2N_{k-1}^{n-1,p} - N_k^{n-1,p}] + N_O(n-1,k,p-1).$$

Substituting the expressions for $N_E(n,k,p)$ and $N_O(n,k,p)$ from Theorems 1 and 2, and rearranging terms, we obtain

$$\binom{n-1}{p}^2 - \binom{n-2}{p}^2 - \binom{n-2}{p-1}^2 = \sum_{i=0}^{k-2} N_i^{n-1,p} - \sum_{i=0}^{k-2} N_i^{n-2,p} - \sum_{i=0}^{k-2} N_i^{n-2,p-1}$$
$$+ \frac{n-1}{k+1} N_k^{n-1,p} + N_{k-1}^{n-1,p} - N_k^{n-1,p}.$$

Applying $\Delta_k$ to both sides of the equation, simplifying, and replacing $n$ by $n+1$ and $k$ by $k-1$ results in

$$\frac{n-k-1}{k+1} N_k^{n,p} = \frac{n-2k}{k} N_{k-1}^{n,p} + N_{k-2}^{n-1,p} + N_{k-2}^{n-1,p-1}, \quad k \geq 1.$$

**Proof of Theorem 5**

From (2) we have for $j \geq 0$

$$\sum_{\substack{n \geq j \\ 0 \leq r \leq n-j}} M_{r,r+j}^{n,k} x^r y^{n-r} = [x+y+2f(x,y)]^k \left[\frac{f(x,y)}{x}\right]^j$$

$$= \sum_{i=0}^{k} \binom{k}{i} (x+f)^i (y+f)^{k-i} \left(\frac{f}{x}\right)^j \quad (6)$$

where $f(x,y)$ satisfies $f=(x+f)(y+f)$. We will express the right hand side of (6) as a power series in x and y by means of Lagrange inversion, as follows:

Let $f=z(x+f)(y+f)$,
$G(t)=(x+t)(y+t)$, and
$\varphi(t) = (x+t)^a (y+t)^b t^c$ .

We use Lagrange inversion in the form (see [2])

$$\varphi(f) = \sum_n [t^n] \left(1 - t\frac{G'(t)}{G(t)}\right) \varphi(t) G^n(t) z^n .$$

Letting z=1, we have

$$(x+f)^a (y+f)^b f^c = \sum_n [t^n] \left(1 - t\frac{G'(t)}{G(t)}\right) \varphi(t) G^n(t) .$$

Since $1 - t\frac{G'(t)}{G(t)} = \frac{xy - t^2}{(x+t)(y+t)}$ , we have

$$(x+f)^a (y+f)^b f^c = \sum_n [t^n] \frac{xy - t^2}{(x+t)(y+t)} (x+t)^a (y+t)^b t^c ((x+t)(y+t))^n$$

$$= \sum_n [t^n] xy (x+t)^{n+a-1} (y+t)^{n+b-1} t^c - \sum_n [t^n] (x+t)^{n+a-1} (y+t)^{n+b-1} t^{c+2} \quad (7)$$

The coefficient of $x^l y^m$ in the right hand side of (7) is seen to be $\binom{n+a-1}{l-1}\binom{n+b-1}{m-1} - \binom{n+a-1}{l}\binom{n+b-1}{m}$ . This term will arise only when the exponent of t is n, i.e., when

$n+a-1-(l-1)+n+b-1-(m-1) = n+a-1-l+n+b-1-m+(c+2) = n$ .

Hence, $n=l+m-a-b-c$ and we have

$$(x+f)^a(y+f)^b\left(\frac{f}{x}\right)^c = \sum_{\substack{l \geq a+c \\ m \geq b+c}} x^{l-c}y^m\left[\binom{l+m-b-c-1}{l-1}\binom{l+m-a-c-1}{m-1} - \binom{l+m-b-c-1}{l}\binom{l+m-a-c-1}{m}\right] \quad (8)$$

To obtain Theorem 5, we first make the substitutions $a=i$, $b=k-i$, $c=j$, $l=r+j$, and $m=n-r$ in the right hand side of (8) to obtain the coefficient of $x^r y^{n-r}$, obtaining

$$(x+f)^i(y+f)^{k-i}\left(\frac{f}{x}\right)^j = \sum_{n,r}\left[\binom{n-k+i-1}{r+j-1}\binom{n-i-1}{n-r-1} - \binom{n-k+i-1}{r+j}\binom{n-i-1}{n-r}\right]x^r y^{n-r} .$$

Theorem 5 follows immediately by letting j=s-r and substituting the result in the right hand side of (6).

### REFERENCES

1. I.Gessel, W.Goddard, W.Shur, H.S.Wilf, and L.Yen, Counting pairs of lattice paths by intersections, J. Combinatorial Theory, Ser. A, to be published.

2. I.M.Gessel, A combinatorial proof of the multivariable Lagrange inversion formula, J. Combin. Theory Ser. A45 (1987), 178-195.

10/15/95